\documentclass{elsart}
\usepackage{graphicx,color}
\usepackage{amssymb,amsmath}
\usepackage{hyperref}

\DeclareMathOperator{\cotan}{cotan}

\newcommand{\eps}{\varepsilon}
\newcommand{\set}[1]{\left\{#1\right\}}

\begin{document}

\begin{frontmatter}

\title{ Low and high frequency approximations to eigenvibrations of  string with double contrasts}

\author[address1]{Natalia Babych\corauthref{cor1}}
\author[address2]{Yuri Golovaty}
\corauth[cor1]{Corresponding author}
\address[address1]{ University of Bath, Bath BA2 7AY, United Kingdom} 
\address[address2]{ Lviv National University, Lviv 79000, Ukraine} 
\ead{n.babych@bath.ac.uk (Natalia Babych)}
\ead{yu$\_$holovaty@franko.lviv.ua}

\begin{abstract}
We study eigenvibrations for inhomogeneous string consisting of two parts
with strongly contrasting stiffness and mass density.
In this work we treat a critical case for the {\it high frequency approximations}, namely
the case when the order of mass density inhomogeneity
is the same as the order of stiffness inhomogeneity,
with heavier part being softer.
The limit problem for high frequency approximations
depends nonlinearly on the spectral parameter.
The quantization of the spectral semiaxies is applied
in order to get a close approximations of eigenvalues as well as eigenfunctions
for the prime problem under perturbation.
\end{abstract}

\end{frontmatter}

\section{Introduction and problem statement}
Models with high contrasts are widely studied
since their unusual properties
give insight into the behaviour
of new meta- and nanomaterials, including those
which already exist
or are reachable nowadays via modern technologies.
The corresponding mathematical problems
often cause computational difficulties
and require new methods of numerical approximation.
A system under consideration possessing two components
with double high contrasts, both in stiffness and mass density,
expresses two distinguishing cases of the limit eigenvibration behaviour
for each of low and high frequency levels.
The description of such systems should not be restricted to
the construction of classical number-by-number eigenfunction asymptotics,
which are called \emph{low frequency approximations}.
They only ensure close approximations to several eigenfunctions
corresponding to the bottom of the spectrum.
For more precise eigenfunction description in the upper part
of the spectrum the classical approach sets the requirement
for $\eps$ to be negligibly small.
Nevertheless, in actual physical models the parameter
$\eps$, denoting the ratio of inhomogeneity for a certain physical
characteristic, is often \emph{small but fixed}.
Then describing actual vibrating systems, a problem of
adequate approximation to eigenfunctions with large numbers arises.
In order to solve the problem we propose a new asymptotics,
being called \emph{high frequency approximations},
and compare them with the classical ones.
The high frequency approximations quite precisely describe
eigenvibrations for which
low frequency approximations are not precise enough.

\textbf{Methods and results.}
Starting from an operator with a discrete spectrum
a classical spectral analysis
provides the discreteness of low frequency limits
for eigenelements of the system with high contrasts.
Nevertheless,
the standard approach misses a certain important characteristic,
because the completeness of eigenfunction system
is lost in the limit.
Accomplishing the investigation
and filling up the gaps in the limit behaviour description
we construct and justify high frequency approximations
to the eigenfunctions. The
quantization conditions play a vital part
in the asymptotics providing an $\eps$--network on the spectral axis
in the range of approximation. Therefore even the leading
terms in the spectral approximations change along
with $\eps$. Thus we obtain a quite precise approximations
to eigenelements
of the prime problem with a fixed small $\eps$.
Comparing to the previous study of the stiff problems \cite{BabychGolovaty1999},
where the leading terms of high frequency approximations are independent
of $\eps$ and the quantization provides a right choice of the correctors,
in the present problem the quantization conditions,
arising in particular from matching WKB and power series expansions,
come along with the choice of the leading terms.

The preliminary results on the limit behaviour of
the system under consideration have been 
discussed in \cite{BabychGolovaty2007b}.
The question of asymptotic description
of low and high frequency eigenvibrations originates
in work \cite{Panasenko} arising again in [3--6]
for problems with perturbations of the stiffness 
only.
Elastic problems with perturbations of stiffness and mass density,
either with other geometries or at different perturbation rates,
have been studied in [6--12].

\textbf{Problem statement.}
Let a stiff
and relatively light part of the string, which occupies an interval
$(a,0)$, be complemented by a flexible and heavy body part occupying $(0,b)$
with $a<0<b$.
We consider a stiffness coefficient
being $k(x)$ on $(a,0)$ and $\eps \varkappa(x)$ on $(0,b)$,
and  mass density being $\eps r(x)$
on $(a,0)$ and $\rho(x)$ on $(0,b)$,
with all functions being positive and smooth in
$[a,0]$ and $[0,b]$
respectively.
We assume that eigenvibrations of
the string are described by the self-adjoint eigenvalue problem
\begin{align}
\label{1DPerturbedProblemEqA}
    &(k(x)u_\eps')' + \eps\lambda^\eps r(x) u_\eps = 0,
    \quad x\in (a,0),
    \qquad u_\eps(a)=0,
    \\
    \label{1DPerturbedProblemEqB}
    \eps &( \varkappa(x) u_\eps')' +
    	\phantom{\eps} \lambda^\eps \rho(x) u_\eps = 0,
    \quad x\in (0,b),
    \qquad u_\eps(b)=0,
    \\
    \label{1DPerturbedProblemCond}
    &u_\eps(-0)=u_\eps(+0),
    \quad
    k(0)u_\eps'(-0)=\eps \varkappa(0)u_\eps'(+0).
\end{align}
We investigate the question how the eigenvibrations of the
media, namely eigenvalues $\lambda^\varepsilon$ and eigenfunctions
$u_\varepsilon$, change if the parameter $\varepsilon$ tends to $0$.
More precisely, we look for the good approximations
of $\lambda^\varepsilon$ and $u_\varepsilon$
as $\eps \to 0$.

\section{Low frequency approximations}
\label{Low frequency approximations}
It is well-known that for each fixed $\varepsilon>0$ the spectrum of problem
(\ref{1DPerturbedProblemEqA})--(\ref{1DPerturbedProblemCond})
is real and discrete,
consisting of simple eigenvalues that
form a sequence
    $0< \lambda^\varepsilon_1 < \lambda^\varepsilon_2 < \dots < \lambda^\varepsilon_n < \dots \to\infty$
    as $n\to\infty$.
The corresponding eigenfunctions $\{ u_{\varepsilon,n}\}_{n=1}^\infty$
 form a basis in $L^2(a,b)$.
Moreover, for each number $n$ the eigenvalue branch $\lambda^\varepsilon_n$ is a continuous function of
$\varepsilon\in(0,1)$ such that $\lambda^\varepsilon_n \le c_n \varepsilon$ with positive constant $c_n$ independent of $\varepsilon$,
which follows from the mini-max principle
since quadratic forms are continuously depending on $\eps$
\cite{Kato}.

Studying the asymptotic behaviour as $\varepsilon\to0$ of each eigenvalue branch $\lambda^\varepsilon_n$ with fixed number $n$
and corresponding eigenfunctions $u_{\varepsilon,n}$,
we immediately have
the convergence
$\eps^{-1 } \lambda^\eps_n \to \lambda_n$
and $u_{\eps,n} \to U_n$, where
$U_n = 0 $ in $(a,0)$
and $U_n$ in $(0,b)$
coincides with an eigenfunction $u^+_n$ of the limit problem
\eqref{low freq 1} for the eigenvalue $\lambda_n$.

We look for the approximations of eigenvalues
and eigenfunctions in the form
\begin{equation}
\lambda^\eps_n
\sim
\eps \mu_n + \eps^2 \nu_n + \dots
,
\quad
u_{\eps,n}
\sim
u_n(x) + \eps w_n(x) +  \dots ,
\ \ x \in (a,b).
\label{low expansions}
\end{equation}
Constructing standardly the asymptotic expansions we
first define the leading terms, which satisfy the problem
\begin{equation}
\text{ for \ } x \in (a,0):
\ \ \ (k(x) u_n' )' = 0
,
\ \
u_n(a) = 0,
\ \
u_n'(-0) = 0.
\label{low problem u^-_0}
\end{equation}
Hence $u_n \equiv 0$ on $(a,0)$ and
therefore
\begin{equation}
\text{ for \ } x \in (0,b):
\ \ ( \varkappa(x) u_n' )' = - \mu_n \rho(x) u_n ,
\quad
u_n(+0) = u_n(b) = 0.
\label{low freq 1}
\end{equation}
Since we are looking for the eigenfunction approximations,
which are supposed to be different from zero,
the limit
$\mu_n$ has to be an eigenvalue with corresponding
eigenfunction $u_n$ of problem
\eqref{low freq 1}.

Let us fix an eigenvalue $\mu_n$ of
\eqref{low freq 1},
and corresponding eigenfunction $u_n$
such that
$
\int_0^b \rho u_n^2 dx = 1.
$
Then the next terms of
\eqref{low expansions}
satisfy the problem
\begin{equation}
\text{ for \ } x \in (a,0):
\ \ (k(x) w_n' )' = 0,
\ \
w_n(a) = 0,
\ \
(k w_n')(-0) =
(\varkappa u_n')(+0).
\label{low problem u^-_1}
\end{equation}
Therefore on $(a,0)$ we have
$
w_n = (\varkappa u_n')(+0) \int_a^x k^{-1}(t) d t
$,
and
\begin{align}
\text{ for \ } x \in (0,b): \ \
&
( \varkappa(x) w_n' )' + \mu_n \rho(x) w_n
=
- \nu_n \rho(x) u_n,
\label{low problem u^+_1: eq}
\\
&
w_n(b) = 0,
\quad
w_n(+0) = v_n(-0).
\label{low problem u^+_1: bc}
\end{align}
The solvability of
\eqref{low problem u^+_1: eq},
\eqref{low problem u^+_1: bc}
along with \eqref{low problem u^-_1} and
the normalization of $u_n$ 
implies
\begin{equation}
\nu_n =
- (k w_n w_n')(-0)
=
- \int_a^0 k (w_n')^2 dx.
\end{equation}
Then on the interval $(0,b)$ we fix a unique $w_n$ such that
$\int_0^b \rho u_n w_n dx = 0$.

\textbf{Justification of low frequency approximations.}
We use the same letter $f$ both for a function defined on
the interval $(a,b)$ and a vector $(f_-,f_+)$, where $f_-$, $f_+$ are the restrictions of
$f$ to $(a,0)$ and $(0,b)$ respectively.
Let $\mathcal{L}$ be the Hilbert space $L^2_r(a,0)\times L^2_\rho(0,b)$ with the scalar product
$
  (u,v)_{\mathcal{L}}=\int_a^0 r u_- v_- d x +
  	\int_0^b \rho u_+ v_+ d x
$
and norm $\|u\|=(u,u)_{\mathcal{L}}^{1/2}$, where $u=(u_-,u_+)$. Let us introduce the matrix operator
$\mathcal{A}_\eps  $ in $\mathcal{L}$
\begin{equation*}
  \mathcal{A}_\eps  =
  \begin{pmatrix}
  -\frac{1}{\eps r}\frac d {dx}(k\frac d {dx})& 0 \\
    0 & -\frac{\eps }\rho\frac d {dx}(\varkappa\frac d {dx})
  \end{pmatrix}.
\end{equation*}
with the domain
\begin{equation*}
 \begin{split}
 \mathcal{D}(\mathcal{A}_\eps  )
 =
 \set{ u \in \mathcal{L}: \right.\ &u_- \in W^2_2(a,0), \ \
  u_-(a)=0,\ \ u_+\in W^2_2(0,b),
 \ \ u_+(b)=0,\\
&\left. u_-(0)=u_+(0),\ \ k(0)u'_-(0)=\eps \varkappa(0)u'_+(0)}.
 \end{split}
\end{equation*}
The $\mathcal{A}_\eps  $ is a self-adjoint operator with a compact resolvent. The spectrum $\sigma(\mathcal{A}_\eps  )$ is the set
of all eigenvalues of (\ref{1DPerturbedProblemEqA})--(\ref{1DPerturbedProblemCond}).

Let $B$ be a self-adjoint operator in Hilbert space $H$ with
a domain $\mathcal{D}(B)$.
Recall that a pair $(\mu, u)\in \mathbb{R}\times\mathcal{D}(B)$
with $\|u\|_H=1$ is a quasimode of the operator $B$
with an accuracy up to $\sigma>0$ if
$\|(B-\mu I)u\|_H\leq \sigma$.

\begin{lem}\label{LemmaVishik}
Suppose that the spectrum of $B$ is discrete. If $(\mu, u)$ is a
quasimode of $B$ with accuracy to $\sigma$, then interval
$[\mu-\sigma,\mu+\sigma]$ contains an eigenvalue of $B$.
Furthermore, if segment $[\mu-\tau,\mu+\tau]$, $\tau>0$, contains one and
only one eigenvalue $\lambda$ of $B$, then $\|u-v\|_H\leq
2\tau^{-1}\sigma$, where $v$ is an eigenfunction of $B$ for the eigenvalue
$\lambda$, $\|v\|_H=1$. \emph{\cite{VishikLusternik1957}}
\end{lem}

\begin{thm}
\label{low freq thm}
For each $n \in \mathbb{N}$
there exists $C_n > 0$ such that
the $C_n \eps^2$--vicinity of
$\eps \mu_n$
contains exactly one eigenvalue $\lambda^\eps_n$
of problem
(\ref{1DPerturbedProblemEqA})--(\ref{1DPerturbedProblemCond}):
\begin{equation}
| \lambda^\eps_n - \eps \mu_n | \le C_n \eps^2.
\label{low eigenvalue error}
\end{equation}
The corresponding normalized eigenfunction $u_{\eps,n}$
satisfies the estimate \newline
$
 \| u_{\eps,n} - u_n - \eps w_n \|_{L_2(a,b)}
 \le \tilde{C}_n \eps^2,
$
with a certain
$\tilde{C}_n > 0$ independent of $\eps$.
\end{thm}

\emph{Proof.}
We introduce a corrector $\phi_n(x) = a^{-1} w_n'(+0) x (x-a)$ on $(a,0)$ and
$\phi_n(x) = 0$ $(0,b)$ such that
$U^\eps_n = u_n + \eps ( w_n + \phi_n)$ belongs to $\mathcal{D}(\mathcal{A}_\eps)$.
Let
$\Lambda^\eps_n = \eps \mu_n + \eps^2 \nu_n$
and
$\tilde{U}^\eps_n = \tau^\eps_n U^\eps_n$
with $\tau^\eps_n = \| U^\eps_n  \|_{\mathcal{L}}^{-1}$.
By the construction
\begin{align}
	& \| \mathcal{A}_\eps   \tilde{U}^\eps_n -
	\Lambda^\eps_n  \tilde{U}^\eps_n  \|_{\mathcal{L}}^2
	\le
	\nonumber
	\\
	& K_1 \eps^4  (\tau^\eps_n)^2
	|\mu_n + \eps \nu_n|^2 ( u_n'(+0)^2 + w_n'(+0)^2)
	+
	K_2 \eps^6 (\tau^\eps_n)^2 \nu_n^2 \| w_n \|^2_{L^2_\rho(0,b)} ,
	\label{precise estimate}
\end{align}
with positive constants $K_j$ independent of $\eps$ and $n$,
and also
\begin{equation}
| \tau^\eps_n | \le ( 1 - \eps \| w_n + \phi_n\|_{\mathcal{L}})^{-1}
\le 1 + \hat{C}_n \eps
\label{tau}
\end{equation}
for $\eps$ small enough.
Therefore
a pair $\Lambda^\eps_n$ and
$\tilde{U}^\eps_n$
is a quasimode of $\mathcal{A}_\eps$ with the accuracy up to $C_n \eps^2$.
By Lemma~\ref{LemmaVishik},
in $C_n \eps^2$--vicinity of $\Lambda^\eps_n $ there exists
a certain eigenvalue $\lambda^\eps_j$ of
(\ref{1DPerturbedProblemEqA})--(\ref{1DPerturbedProblemCond}).
Additionally, it can be easily shown that
the eigenvalues converge
saving multiplicity, $\eps^{-1}\lambda^\eps_n \to \mu_n$.
Since the limit problem has only simple eigenvalues,
in a certain $\hat{C}_n\eps$--vicinity
of $\mu_n$ there is no other eigenvalues of
(\ref{1DPerturbedProblemEqA})--(\ref{1DPerturbedProblemCond})
that provides \eqref{low eigenvalue error}.
Applying again Lemma~\ref{LemmaVishik}
finishes the proof.
\hfill $\Box$

Note that low frequency vibrations vanish
in $(a,0)$ as $\varepsilon\to 0$. This naturally raises the question on the possibility of constructing
other \emph{non-trivial on} $(a,0)$ approximations of eigenvibrations addressed next.

\section{High frequency approximations}
\label{Sec4}
Considering sufficiently large eigenvalues
$
\lambda^\eps_n\sim \eps^{-1}(\omega+\eps\omega_1)^2
$
with $\omega>0$,
we look for the asymptotic expansions of  eigenfunctions
$u_{\eps,n}(x) \sim Y(\eps,x)$ with
\begin{align}
Y(\eps,x)
=
\left\{
\begin{array}{ll}
v_0(x)+\eps v_1(x)+\eps^2 v_2(x),
&
\ x \in (a,0),
\\
      \left(c_0(x)+\eps c_1(x)\right)
      \sin \gamma_\eps S(x)
      +
      \eps c_2(x)\cos \gamma_\eps S(x),
  		&
\ x \in (0,b),
\end{array}
\right.
\label{1DHighFrequencyAppr}
\end{align}
where $v_0$ is different from zero and
$
\gamma_\eps = \frac{\omega}{\eps}+\omega_1.
$
The expansion in form \eqref{1DHighFrequencyAppr}
consists of power series on the interval $(a,0)$
and two-term short-wave (WKB) approximation \cite{Fedoryuk}
on $(0,b)$ since equation \eqref{1DPerturbedProblemEqB}
contains a small parameter near the
highest derivative.
Substituting these expressions into
equation and boundary condition
\eqref{1DPerturbedProblemEqA} gives
\begin{align}
		\label{1DproblemV0}
    &(kv_0')'+\omega^2 rv_0=0, \qquad v_0(a)=0,
    \\
    \label{1DproblemV1}
    &(kv_1')'+\omega^2 rv_1=-2\omega\omega_1 rv_0,\qquad v_1(a)=0,
    \\
    \label{1DproblemV2}
    &(kv_2')'+\omega^2 rv_2=-\omega_1^2 rv_0-2\omega\omega_1 rv_1,
    \qquad v_1(a)=0.
\end{align}
Next, we substitute $Y(\eps,x)$ into \eqref{1DPerturbedProblemEqB}:
\begin{align}
   \eps\gamma_\eps^2
   &
   \left(-\varkappa {S'}^2
   +\rho\right) Y(\eps,\cdot)
   +
   \eps\gamma_\eps\Bigl(2\varkappa S'c_0'+(\varkappa S')'c_0\Bigr)
   \cos \gamma_\eps S
   \nonumber
   \\
   &+
   \eps\omega\Bigl(2\varkappa S'c_1'+(\varkappa S')'c_1\Bigr)\cos \gamma_\eps S
   \nonumber
   \\
   &-
    \eps\Bigl(2\omega\varkappa S'c_2'
    +\omega(\varkappa S')'c_2-(\varkappa c_0')'\Bigr)
    \cos \gamma_\eps S
   =O(\eps^2),
\label{1DSubstitutionOfWKB}
\end{align}
Equating the expressions in the large parentheses to zero
we minimize the discrepancy in \eqref{1DSubstitutionOfWKB}.
The \emph{eikonal} equation $\varkappa {S'^2}=\rho$ has a solution
\begin{equation*}
    S(x)=\int^b_x \varkappa^{-1/2}(\tau)\rho^{1/2}(\tau)\,d\tau,\quad x\in (0,b).
\end{equation*}
Consequently, the \emph{transport} equation
$2\varkappa S'c'+(\varkappa S')'c=0$
admits a solution
$c(x)=\varkappa^{-1/4}(x)\rho^{-1/4}(x)$
up to a constant multiplier.
Therefore $c_0(x)=\beta_0 c(x)$ and $c_1(x)=\beta_1 c(x)$.
Introducing $h$ as a unique solution of the 
problem
$$
2\varkappa S'h'+(\varkappa S')'h=(\varkappa c')'
\ \text{for} \ x < b,
\quad h(b)=0,
$$
we set $c_2(x)=\beta_0\omega^{-1}h(x)$
providing the boundary condition $Y(\eps,b)=0$ is satisfied.
By construction
$Y(\eps,\cdot)$
formally solves  equation \eqref{1DPerturbedProblemEqA} up to
the terms of order $\eps^3$ and equation \eqref{1DPerturbedProblemEqB} up to the terms of order $\eps^2$.

We now apply interface conditions \eqref{1DPerturbedProblemCond}
in order to define parameters
$\omega$, $\omega_1$, $\beta_0$ and $\beta_1$.
Before that, regularizing the $\eps$-dependence of $Y(\eps,+0)$ we apply
the restriction
\begin{equation}
\label{1DPreQuantizingCondition}
    \left(\frac{\omega }{\eps}+\omega_1\right)S(0)= \delta+\pi l,
    \qquad \delta\in (-\pi/2,\pi/2],
    \quad l\in \mathbb{Z} .
\end{equation}
Satisfying the interface conditions up to the terms of order $\eps^2$,
we set
\begin{align}
&    \begin{cases}
      \phantom{k(0)} v_0(0)=(-1)^l\beta_0 c(0)\sin\delta\\
      k(0)v'_0(0)= (-1)^l\beta_0 \omega S'(0)\varkappa(0) c(0)\cos\delta,
    \end{cases}
\label{1DSystemsV0(0)}
    \\
&    \begin{cases}
        \phantom{k(0)}v_1(0)=(-1)^l(\beta_1 c(0)\sin\delta+g_1)\\
      k(0)v'_1(0)= (-1)^l(\beta_1\omega S'(0)\varkappa(0)c(0)\cos\delta +\beta_0\varkappa(0)g_2),
    \end{cases}
\label{1DSystemsV1(0)}
\end{align}
where $g_1=\beta_0\omega^{-1}h(0)\cos\delta$,
$g_2=\omega_1S'(0)c(0)\cos\delta+(c'(0)-S'(0)h(0))\sin\delta$.

Combining \eqref{1DproblemV0} and \eqref{1DSystemsV0(0)}
 we obtain that $v_0$
is a solution to the problem
 \begin{align}
 \label{1DPencilV0}
 \begin{array}{c}
    (kv')'+\omega^2 rv=0
    \ \ \text{in} \ \ (a,0), \\
    v(a)=0,\quad k(0)v'(0)\sin\delta-\omega \varkappa(0)S'(0)v(0)\cos\delta=0.
    \end{array}
 \end{align}

\begin{prop}
\label{1DdeltaExistence}
  For every $\omega>0$ there exists a unique  $\delta(\omega)\in (-\pi/2,\pi/2]$ such that problem \eqref{1DPencilV0} has a nontrivial solution $v$.
\end{prop}

\emph{Proof.}
If $\omega^2$ is an eigenvalue of the problem $(kv')'+\omega^2 rv=0$, $v(a)=0$, $v(0)=0$, we put $\delta(\omega)=0$.
Otherwise we consider the eigenvalue problem
\begin{equation}
\label{1DauxiliaryProblem}
    (k v')'+\omega^2 r v=0
     \ \ \text{in} \ \ (a,0),
     \quad v(a)=0,\quad k(0)v'(0)+\mu v(0)=0
\end{equation}
with respect to the spectral parameter $\mu$. For each $\omega$  under consideration the problem has a unique eigenvalue
$\mu(\omega)$, which is due to the fact that the spectral parameter is missed in equation. Therefore $\delta(\omega)$ can be found as a unique root in $(-\pi/2,\pi/2]$ of the equation
\begin{align}
\omega \varkappa(0)|S'(0)| \cotan \delta=\mu(\omega).
\label{mu(omega)}
\end{align}
Recall that $S'(0)<0$.
\hfill$\Box$

Fixing an arbitrary $\omega > 0$ we also fix $v_0=v(\omega, x)$
being corresponding eigenfunction of non-linear pencil
\eqref{1DPencilV0} with $\delta = \delta(\omega)$ defined by
Proposition~\ref{1DdeltaExistence}.
Let additionally $v_0$ be unity normalized in $L^2_r(a,0)$.
Consequently,
\eqref{1DSystemsV0(0)} provides
$$
\beta_0=
\left\{
\begin{array}{ll}
\frac{(-1)^l v_0(0)}{c(0)\sin\delta(\omega)}
& \text{if} \ \ \delta(\omega)\neq 0,
\\
\frac{(-1)^l k(0) v'_0(0)}{\omega \varkappa(0)S'(0) c(0)}
& \text{if} \ \ \delta(\omega)=0.
\end{array}
\right.
$$
We conclude from condition \eqref{1DPencilV0} that the function $\beta_0(\omega)$ is continuous at every point $\omega_*$ for which $\delta(\omega_*)=0$.
 From \eqref{1DSystemsV1(0)} and
\eqref{mu(omega)} we obtain
 \begin{align}
 \label{1DPencilV1}
 \begin{array}{ll}
    (k v_1')'+\omega^2 r v_1=-2\omega\omega_1r v_0
    \ \ \text{in} \ \ (a,0),
    \\
    v_1(a)=0,\quad k(0)v'_1(0)+\mu(\omega)v_1(0)=\varkappa(0)f,
    \end{array}
 \end{align}
where $f=g_2\sin\delta(\omega)-\omega S'(0)g_1\cos\delta(\omega)$ and $\delta(\omega)\neq 0$. The problem admits a solution if and only if $v_0(0)\varkappa(0)f=-2\omega\omega_1$, since $\mu(\omega)$ is an eigenvalue of \eqref{1DauxiliaryProblem}.
This solvability condition can be derived multiplying the equation by $v_0$ and integrating twice by parts.
It may be written in the form
\begin{equation*}
    \omega_1=
    \frac{\left(h(0)S'(0) - c'(0) \sin^2\delta(\omega)\right)v^2(\omega,0)}
  {\left(2\omega + \varkappa(0)S'(0)\cos\delta(\omega)\right)c(0)\sin\delta(\omega)}
  \ \ \ \text{if} \ \ \ \delta(\omega) \neq 0.
\end{equation*}
Thus we get $\omega_1$ as a function of $\omega$.
Additionally,
we  obtain
$$
\omega_1=-k(0)v_0'(0)c_2(0) (2\omega)^{-1}
\ \ \ \text{if} \ \ \ \delta(\omega)=0.
$$
We now can find $v_1$, which is ambiguously determined.
Subordinating it to the condition
$ \int_a^0 r v_0 v_1 \, d x = 0$ we fix it uniquely.
Then $\beta_1$
is given by \eqref{1DSystemsV1(0)}.
We fix an arbitrary $v_2$ being solution of \eqref{1DproblemV2}.

Let us  return to condition \eqref{1DPreQuantizingCondition}.
Now it may be considered as the countable set of equations for $\omega$:
\begin{equation}
\label{1DQuantizingCondition}
    \left(\frac{\omega }{\eps}+\omega_1(\omega)\right)S(0)-\delta(\omega)=\pi l, \qquad l\in \mathbb{Z}.
\end{equation}
Since $\omega_1(\omega)$ can have a vertical asymptote
in the interval $I = [0,\frac12 \varkappa(0)|S'(0)|)$,
equation \eqref{1DQuantizingCondition} can have
roots in $I$.
More subtle analysis shows that for each $l$ there always exists a unique root of \eqref{1DQuantizingCondition} in the set
$[\frac12 \varkappa(0)|S'(0)|,\infty)$ because
$\omega_1(\omega)\to 0$ as
$\omega\to +\infty$ and quantity $\delta(\omega)$ is bounded.
We consider the roots that increase along with $l$.

\begin{defn}
  We say that $\omega(l)$ is an \emph{admissible limit frequency} for given $\eps>0$ and $l\in \mathbb{Z}$ if it is the largest root of \eqref{1DQuantizingCondition}.
\end{defn}

Let us establish connection between the exact eigenfrequencies
$\sqrt{\lambda^\eps_{l}}$
and addmissible limit frequencies  $\omega(l)$ in case
of constant coefficients.
Indeed, for $k=\varkappa=r=\varrho=1$ we have
$\sqrt{\lambda^\eps_{n}} = \sqrt{\eps} \pi n (b-\eps a)^{-1}$,
$ n\in \mathbb{N} $.
Counting the admissible frequencies in this case we note that
$S(x) = b-x$ and $v_0(x) = C_0 \sin \omega (x-a)$
providing, via the proof of Proposition \ref{1DdeltaExistence},
$\delta(\omega) = 0$ if $\omega = \pi n a^{-1}$ for natural $n$
and $\cotan \delta(\omega) = \omega^{-1} \mu(\omega)$ for all other $\omega$. Moreover, \eqref{1DauxiliaryProblem} yields
$\mu(\omega) = -v_0'(0) v_0(0)^{-1}=\omega \cotan (\omega a)$
gaining $\delta(\omega) = \arctan (\tan (\omega a))
\in (-\pi/2, \pi/2)$ or $\delta = \pi /2$.
Observe that $c_0(x)=\beta_0$,
$h(x)=0$ and $c_2(x)=0$ providing $\omega_1 = 0$.
Then \eqref{1DQuantizingCondition} becomes
${\omega } b {\eps}^{-1} = \arctan (\tan (\omega a)) + \pi l $.
Therefore,
${\omega } b {\eps}^{-1} \in (-\pi/2 + \pi l, \pi/2 + \pi l]$
and
$\tan ({\omega } b {\eps}^{-1} ) = \tan (\omega a)$
providing
$
\omega = \eps \pi k(l) (b-\eps a)^{-1}
$
for $k(l) \in \mathbb{N} \cap K^\eps_l$ with
$K^\eps_l = ( z_\eps(l-1/2), z_\eps(l+1/2)]$,
where
$z_\eps = (1 + \eps a (b-\eps a)^{-1})^{-1}$.
Since $z_\eps > 1$ and therefore the length $|K^\eps_l|$
is also larger then $1$, we have at least one natural
$k(l) \in K^\eps_l$.
Picking up the maximal value
$k^{max}_l(\eps) \in K^\eps_l \cap \mathbb{N}$ we fix
the admissible frequency
$\omega(l) =  \eps \pi k^{max}_l(\eps) (b-\eps a)^{-1}$.
Note that $k^{max}_l(\eps) = l $ for the range of numbers
$l < \frac{b+\eps a}{2 \eps |a| }$.
Therefore, $\sqrt{\lambda^\eps_l} =  \frac{\omega(l)}{\sqrt{\eps}}$ for
$l < \frac{b+\eps a}{2 \eps |a| }$.

Having exact correspondence for the range of
eigenfrequencies and admissible frequencies in case
of constant coefficients, in general case
we further use the set
of admissible frequencies as the first approximation
for the eigenfrequencies.
Let $\Phi_\eps$ denote the set of all admissible limit
frequencies.
 The subset $\Phi_\eps$ of $\mathbb{R}_+$ is thick enough, the distance between neighboring roots is  comparable
with $\eps$. In some sense \eqref{1DQuantizingCondition} could be regarded as a kind of WKB quantization condition. The positive spectral ray $\omega>0$ is covered by the $\eps$-net $\Phi_\eps$, for each point of which we can construct the asymptotics \eqref{1DHighFrequencyAppr}. For each admissible frequency $\omega\in \Phi_\eps$ we will denote by
$Y_\omega(\eps,x)$ the corresponding asymptotic solution \eqref{1DHighFrequencyAppr}.

\section{Justification of high frequency approximations}
\label{Sec5}
The function $Y_\omega(\eps,x)$ can be used to construct a quasimode of the operator $\mathcal{A}_\eps $.
Clearly, $Y_\omega\in \mathcal{L}$, but
$Y_\omega\not\in \mathcal{D}(\mathcal{A}_\eps  )$
because of discontinuity at $x=0$.
Let us introduce functions $\zeta_0, \,\zeta_1\in C^1(a,0)$ such that $\zeta_0(a)=0$, $\zeta_0(0)=1$ and
$\zeta_1(a)=0$, $\zeta_1(0)=0$, $\zeta'_1(0)=1$.
Both functions are extended by zero into $(0,b)$.
Introducing
$$
\tau_0(\eps)
=
\left(Y_\omega(\eps,+0)-Y_\omega(\eps,-0)\right)\eps^{-2},
\ \
\tau_1(\eps)
=
\left(\eps Y'_\omega(\eps,+0)-Y'_\omega(\eps,-0)\right)\eps^{-2},
$$
which are bounded in $\eps$ by construction,
we obtain that the function
$$
\tilde{Y}_\omega(\eps,\cdot)
=
Y_\omega(\eps,\cdot)
+
\eps^2(\tau_0(\eps)\zeta_0
+
\tau_1(\eps)\zeta_1)
$$
belongs to $\mathcal{D}(\mathcal{A}_\eps)$.
Setting
$
\Upsilon_\omega(\eps,\cdot)
=
\|\tilde{Y}_\omega\|_{\mathcal{L}}^{-1}\cdot\tilde{Y}_\omega(\eps,\cdot)
$
we prove
the following estimate
$\left\|(\mathcal{A}_\eps -\eps\gamma^2_\eps  I)\Upsilon_\omega(\eps,\cdot)\right\|_{\mathcal{L}}\leq C\eps^2$
by recalling $\eps^{-1}(\omega+\eps \omega_1)^2=\eps \gamma^2_\eps$.

\begin{prop}\label{1DPropQuasimode}
The pair $(\eps\gamma^2_\eps,\Upsilon_\omega(\eps,\cdot))$
is a quasimode of $\mathcal{A}_\eps$ with  accuracy to $O(\eps^2)$
for every admissible frequency $\omega\in \Phi_\eps$.
\end{prop}

\begin{prop}
\label{lemma: number growth}
For the range of numbers $n \le \theta \eps^{\sigma - 1/2}$
with arbitrary $\theta > 0$ and $0< \sigma < 1/2$,
the eigenvalues satisfy the estimate
$ | \lambda^\eps_n - \eps \mu_n| \le K_* \eps^{ 1+2 \sigma }$.
\end{prop}

\emph{Proof.}
In order to improve \eqref{low eigenvalue error}
we calibrate \eqref{precise estimate}.
Eigenfrequencis
$\eta_n = \mu_n^{1/2}$ and normalized eigenfunctions $u_n$
of problem~\eqref{low freq 1}
cam be represented as \cite{Fedoryuk}
\begin{align}
\label{1DAsymptoticsMuN}
    \eta_n & = {\textstyle
    \frac{\pi n}{S(0)} +\frac{\pi S(0)}{n} +
    O(\frac{1}{n^3})}
    \text{ \ \ \ as \ \ } n\to +\infty,
    \\
    u_n
    & =
    \sqrt{2 / S(0)} (\varkappa \rho )^{-1/4}(\sin \eta_n S)(1 + O(\eta_n^{-1}))
    \text{ \ \ on \ } (0,b),
    \label{u n WKB}
\end{align}
where \eqref{u n WKB}  is uniform on $[0,b]$
and admits differentiation in $x$.
Then we have the approximation of the right-hand side in
\eqref{low problem u^+_1: eq}, \eqref{low problem u^+_1: bc}
\begin{equation}
	v_n(-0) = \eta_n \beta_1 \beta_2 (1 + O(\eta_n^{-1})),
	\quad
	\nu_n = -\eta_n^2 \beta_1^2 \beta_2 ( 1 + O(\eta_n^{-1})),
	\quad n \to \infty,
	\label{nu 1}
\end{equation}
with $\beta_1 = \sqrt{2 / S(0)} (\varkappa \rho )^{-1/4}(0) $ and
$\beta_2 = \int_a^0 k(t)^{-1} d t $.
Since there exists the fundamental set of solutions corresponding
\eqref{low problem u^+_1: eq} in the form
 \cite{Fedoryuk}
$$
y_n = ( 1 + O(\eta_n^{-1}))(\varkappa \rho )^{-1/4}\sin \eta_n S
\text{ \ \ and \ \ }
g_n = ( 1 + O(\eta_n^{-1}))(\varkappa \rho )^{-1/4}\cos \eta_n S,
$$
$w_n$ admits representation
$w_n = p_n(x) y_n(x) + q_n(x) g_n(x)$ for certain functions $p_n$ and $q_n$. Exploring this structure of solution in problem
\eqref{low problem u^+_1: eq}, \eqref{low problem u^+_1: bc}
we obtain
\begin{equation}
 w_n(x) = ( 1 + O(\eta_n^{-1}))
 \eta_n K (\varkappa \rho )^{-1/4}(x)\cos \eta_n S(x),
 \ \ \ n \to \infty,
\label{w n}
\end{equation}
with constant $K = \beta_1 \beta_2 (\varkappa \rho)^{1/4}(0) S(0)^{-1}$.
Then \eqref{w n} and \eqref{1DAsymptoticsMuN} provide
\begin{equation}
	\| w_n \|_{L^2_\rho(0,b)} \le K_1 n,
	\ \ \
	| w_n'(+0) | \le K_2 n.
	\label{estimate w n}
\end{equation}
Finally, counting \eqref{1DAsymptoticsMuN}--\eqref{estimate w n}
in \eqref{precise estimate} we obtain
\begin{align}
	&
	\| \mathcal{A}_\eps   \tilde{U}^\eps_n -
	\Lambda^\eps_n  \tilde{U}^\eps_n  \|_{\mathcal{L}}^2
	\le
	K_3 |\tau^\eps_n|	\eps^2 n^2
	\sqrt{ 1 +  \eps^6 n^6 }
	\le
	K_4 \eps^{1+2 \sigma},
	\label{precise 2}
\end{align}
for $n \le \theta \eps^{\sigma -1/2}$ and
$ |\tau^\eps_n| \le 1 - K_5 \eps^{\sigma + 1/2}$, which follows from \eqref{tau}.
Then application of Lemma \ref{LemmaVishik} finishes the proof.
\hfill $\Box$

\begin{thm}
\label{theorem: HFV justification}
Let $ \theta \geq 1 $, $ 0 < \sigma < 1/2 $,
$ 0 < \gamma < 1/2 - \sigma $.
If $\omega = \omega(n)$ is an admissible limit frequency
from the number range
$n \in [\theta^{-1}\eps^{-\gamma},
\theta \eps^{\sigma-1/2}]$
and $\delta(\omega) \not = \pi/2$
then the eigenvalue $\lambda_n^\eps$
and eigenfunction $y_{\eps,n}$ satisfy the estimates
\begin{equation*}
    |\lambda_n^\eps - \eps ^{-1}(\omega+\eps\omega_1(\omega))^2|
    \leq \alpha_1\eps^2,
    \qquad
     \|y_{\eps,n}-\Upsilon_\omega(\eps,\cdot)\|_{L_2(a,b)}
     \leq \alpha_2 \theta\eps^{1+\gamma},
\end{equation*}
with positive constants $\alpha_1$, $\alpha_2$ being independent of $\eps$.
\end{thm}

\noindent
\emph{Proof.}
Let $\omega_n^\eps = \sqrt{\lambda_n^\eps}$.
Proposition \ref{lemma: number growth}
and  \eqref{1DAsymptoticsMuN}
for the given number range yield
\begin{equation}
\label{1DEigenvalueAsymptotics}
\textstyle
 \omega_n^\eps
 =
 \eps^{1/2} \pi \left(\frac{n }{S(0)}
 +\frac{ S(0)}{n}\right)
 +
 O(\eps^{\frac12 + 2 \sigma +\gamma})
 \text{ \ \ and \ \ }
 \lambda_n^\eps =
 \frac{\eps \pi^2n^2}{S^{2}(0)} + O(\eps).
\end{equation}
We now estimate
the distance between neighboring eigenvalues of $\mathcal{A}_\eps$.
From \eqref{1DEigenvalueAsymptotics} we have
$
\lambda_{n+1}^\eps - \lambda_n^\eps =
\eps \pi^2(2n+1)S^{-2}(0) + O(\eps)$. If $n\geq \theta^{-1}\eps^{-\gamma}
$
then
\begin{equation}
\label{1DDistanceBetweenLambda}
    | \lambda_{n+1}^\eps - \lambda_n^\eps |
    \geq
    2 \pi^2 S^{-2}(0) n \eps + O(\eps)
    \geq
    g_0 \theta^{-1}\eps^{1-\gamma},
\end{equation}
with constant $g_0$
being
positive and independent of $n$.
By the similar argument,
\begin{equation}
\label{1DDistanceBetweenOmega}
\textstyle
    | \omega_{n+1}^\eps - \omega_n^\eps |
    \geq
    \sqrt{\eps}(\frac{\pi}{S(0)} - c( n^{-2} +
    \eps^{2 \sigma + \gamma}))
    \geq \sqrt{\eps}\left(\frac{\pi}{S(0)}-\theta_1\eps^{\gamma}\right).
\end{equation}
Let for a certain number $l$ the admissible frequency
$\omega_* = \omega(l)$ minimize the difference
$
\left| \sqrt{\eps} \omega_n^\eps -
( \omega_* + \eps\omega_1(\omega_*) )
\right|
$,
which
equals
$| \eps \pi ( n-l - \frac{\delta}{\pi}) + O(\eps^{1+\gamma})|$
by \eqref{1DQuantizingCondition} and
\eqref{1DEigenvalueAsymptotics}.
Note that if $|\delta| < \frac{\pi}{2}$
the latter is minimized only for $l=n$
(for sufficiently small $\eps$)
and then we have
\begin{align}
    &\left|\sqrt{\eps}
    \omega_n^\eps -
    (\omega_*+\eps\omega_1(\omega_*))\right|
    \le
    |\delta(\omega_*)| S^{-1}(0)\eps + \theta_2 \eps^{1+\gamma}.
\label{1DMinimization}
\end{align}
Suppose that frequency $\omega_{n+1}^\eps$ (the same for $\omega_{n-1}^\eps$) also
sa\-tis\-fies the last inequality.
Then we obtain the estimate
$$
| \omega_{n+1}^\eps - \omega_n^\eps |
\leq
2 | \delta(\omega) | S^{-1}(0) \sqrt{\eps} +
2 \theta_2\eps^{1+\gamma}
$$
that contradicts \eqref{1DDistanceBetweenOmega} for  $\eps<\eps_0$
because $2| \delta(\omega) | < \pi$.
The number $\eps_0$ can be found from the equation
$ S(0) \eps^\gamma (2 \theta_2 \eps^{1/2} + \theta_1 )
= \pi - 2 |\delta(\omega)| $.
Hence, $\omega_n^\eps$ is a unique   eigenfrequency that satisfies
\eqref{1DMinimization} with $\omega(n)$ being the root of
\eqref{1DQuantizingCondition}.

In view of Lemma \ref{LemmaVishik} and Proposition \ref{1DPropQuasimode}, we can improve inequality \eqref{1DMinimization} to
$
|\lambda_n^\eps-\eps ^{-1}(\omega+\eps\omega_1(\omega))^2|\leq \alpha_1\eps^2.
$
Repeated application of Lemma  \ref{LemmaVishik} enables us to write
$
\| u_{\eps,n} - \Upsilon_\omega(\eps,\cdot) \|_{\mathcal{L}}
\leq \alpha_2 \theta\eps^{1+\gamma},
$
because the spectral gap is of order
$\theta^{-1}\eps^{1-\gamma}$ due to \eqref{1DDistanceBetweenLambda}.
\hfill
$\Box$

Note that the case $\delta = \pi / 2$ is not a typical situation.
Indeed, in this case the admissible frequency $\omega$ coincides
with the eigenfrequency of problem \eqref{1DPencilV0}
with Neumann condition $v'(0) = 0$.
In this case we can not establish the vicinity
of $\omega$ that would contain one and only one
eigenfrequency $\sqrt{\lambda^\eps_n}$. The arguments of the last proof
show only that there exist not more then two eigenfrequencies
satisfying \eqref{1DMinimization}.
If $\omega_n^\eps$ and $\omega_{n+1}^\eps$ satisfy  \eqref{1DMinimization}
then $\Upsilon_\omega(\eps,\cdot)$ is not yet a good
approximation to any of $u_{\eps,n}$ or
 $u_{\eps,n+1}$. The situation could be improved by the next terms of asymptotics but that is out of scope of this paper.


\textbf{Numerical example.}
Let us consider coefficients
$k = 1$ and $r = 1 + x^2$ on the interval $(-1,0)$,
$\varkappa = 1$ and $\rho = 1 + x$ on $(0,1)$.
For the value of small parameter
$\eps=0.05$, in Fig. 1 
we have plotted  the eigenfunctions  of \eqref{1DPerturbedProblemEqA}-\eqref{1DPerturbedProblemCond} and the leading terms of the low and high frequency approximations given by \eqref{low expansions} and \eqref{1DHighFrequencyAppr}.
Let us emphasize that the purpose of high frequency approach
is in a good approximation of eigenfunctions.
In the example only the first few eigenfunctions
could be well approximated by the low frequency approach.
Already $u_{\eps,5}$ is quite far away from its
low frequency limit (see Fig. 1), 
and that can not be improved by the next terms of low frequency asymptotics, because the absolute error  is large enough.
Note that the eigenvalue $\lambda^\eps_{5}$ is still far from zero and
thus $\sqrt{\lambda^\eps_{5}}$ can not be treated as a low frequency. In  the right-hand side plots we observe that the high frequency approximations
work well for the range of numbers between $5$ and $15$.
We have to mention that the proof of Theorem \ref{theorem: HFV justification}
is done by asymptotic methods, so it would be challenge
to tell in particular examples the exact range of numbers $n$, for
which high frequency approximations are valid in case
of fixed $\eps$. 

We refer to the values of  $\sqrt{\lambda_n^\eps}$
that are calculated with high accuracy as to ``exact".
The numerical values under discussion
are represented in the Table
\begin{center}
  \begin{tabular}{|c|c|c|c|}
\hline
 { $n$}
 &   \textit{5} & \textit{10} & \textit{15}
\\
\hline
  \text{exact} \ \ $\sqrt{\lambda_n^\eps}$
  &  2.76675   & 5.52678   &  8.27450
\\[-5pt]
  \text{low freq. approximation} \ \ $\sqrt{\eps \mu_n}$
  &  2.88055  & 5.76252   &  8.64418
\\[-5pt]
  {$\omega $}        &  0.6270   & 1.260    &  1.860
\\[-5pt]
   {$\omega_1$}     & -0.22224  & -0.53779 &  -0.02669
\\[-5pt]
   {$\delta$  }     &  -0.63509 & -1.3217  & -0.03770
\\[-5pt]
  \text{high freq. approximation}
  \ \ $\frac{\omega}{\sqrt{\eps}} + \sqrt{\eps}\omega_1 $
  & 2.75433	&	5.51464 & 8.3122
\\
\hline
   \end{tabular}
\\[15pt]
\end{center}
As for numerical example we present the low and high frequency approximations to
eigenfunctions by the leading terms of the expansions only.
Thus the low frequency approximations $\sqrt{\eps \mu_n}$
to eigenfrequencies $\sqrt{\lambda_n^\eps}$
are given in the Table and are accomplished
by visualization of eigenfunctions $u_n$ of problem \eqref{low freq 1}
in the left columns of Fig. 1 
($u_n$ is extended by zero to $(-1,0)$).
In order to find the admissible frequency $\omega$ we need also
$\delta$ and $\omega_1$. In the vicinity of expected $\omega$
we create network over $\omega$ and for each of this
$\omega$ we find $\delta(\omega)$
satisfying \eqref{1DPencilV0}
(up to $10^{-7}$) and then find $\omega_1(\omega)$
such that \eqref{1DPencilV1} has a solution.
Finally, we
find the admissible frequency $\omega=\omega(n)$
giving the best approach to \eqref{1DQuantizingCondition}
over tabulated $\omega$.
The high frequency approximations to the eigenfunctions
we depict from  $u_{\eps,l} \sim v(\omega, x)$
for $x\in(-1,0)$ and
$u_{\eps,l} \sim c_0(x)\sin (\frac{\omega}{\eps}+\omega_1) S(x)$
for $x\in(0,1)$ with
$v(\omega, x)$ from \eqref{1DPencilV0},
$
c_0(x)=\frac{v(\omega,0)}{\sin \delta(\omega)\sqrt[4]{1+x}}
$
and
$ S(x)= 
\frac23(2\sqrt{2}-\sqrt{(1+x)^3})
$.
All depicted eigenfunctions $u_{\eps,n}$ are normalized in $L^2(-1,1)$.

  \begin{figure}[t]
  \begin{center}
   {\scriptsize Low frequency approximations}%
   \hglue70pt{\scriptsize High frequency approximations}\\[5pt]
   \includegraphics[width=0.4\textwidth]{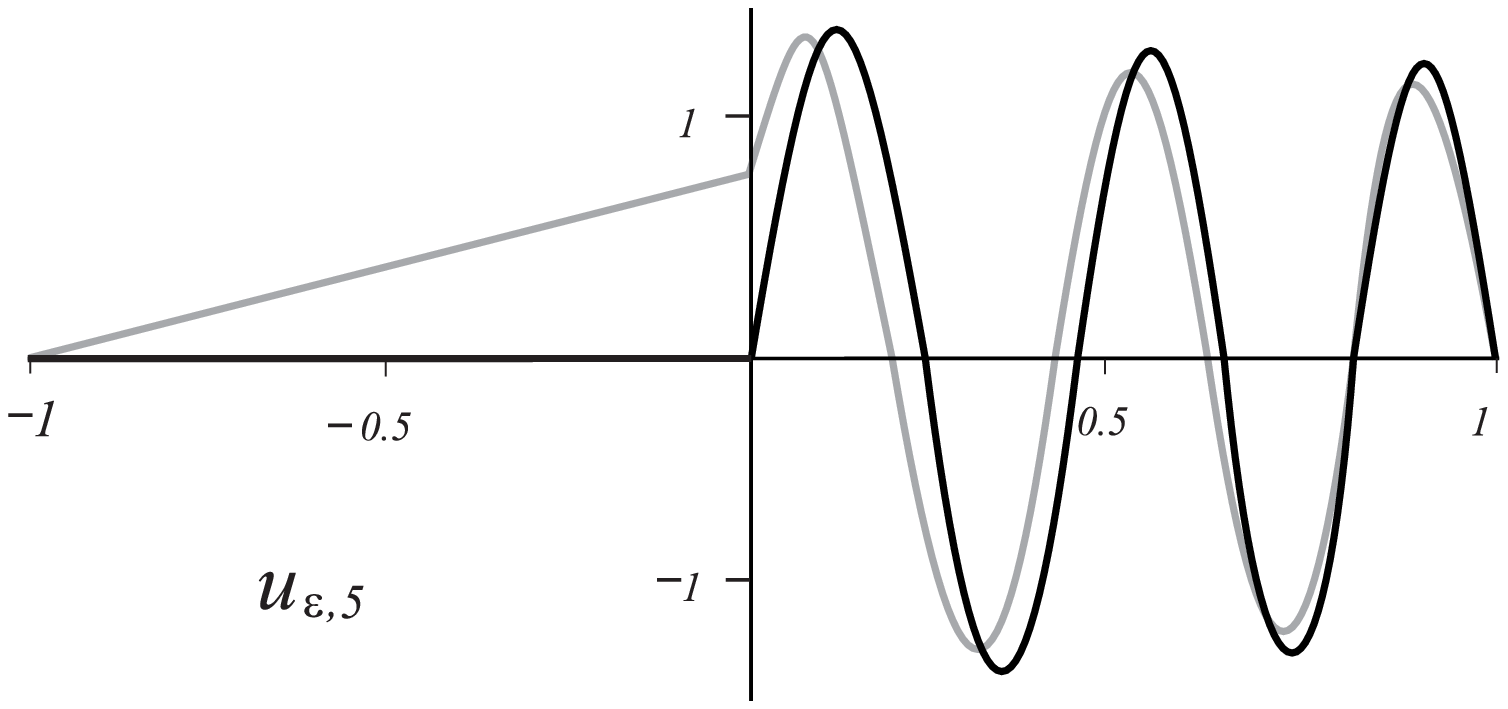}%
   \hglue40pt \includegraphics[width=0.4\textwidth]{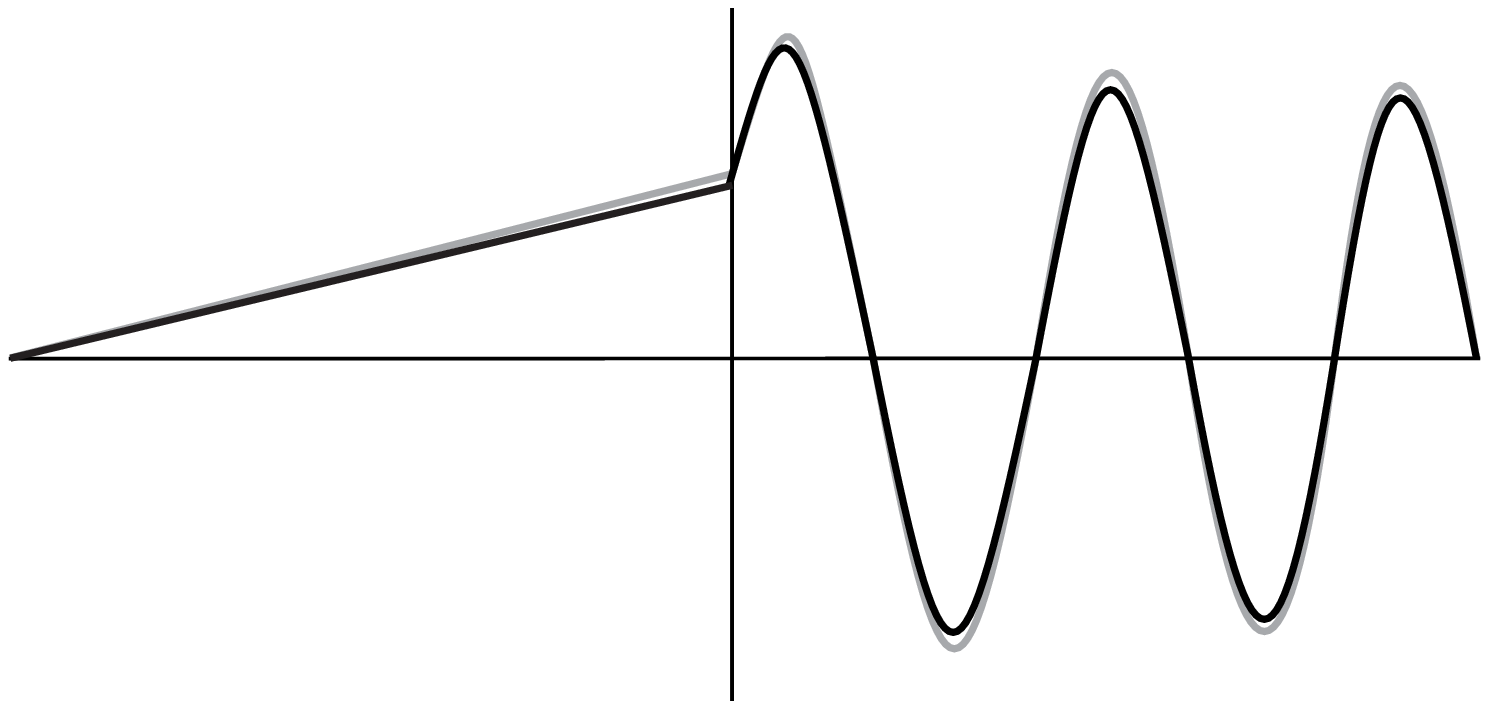}
   \\[-5pt]
\vglue12pt
   \includegraphics[width=0.4\textwidth]{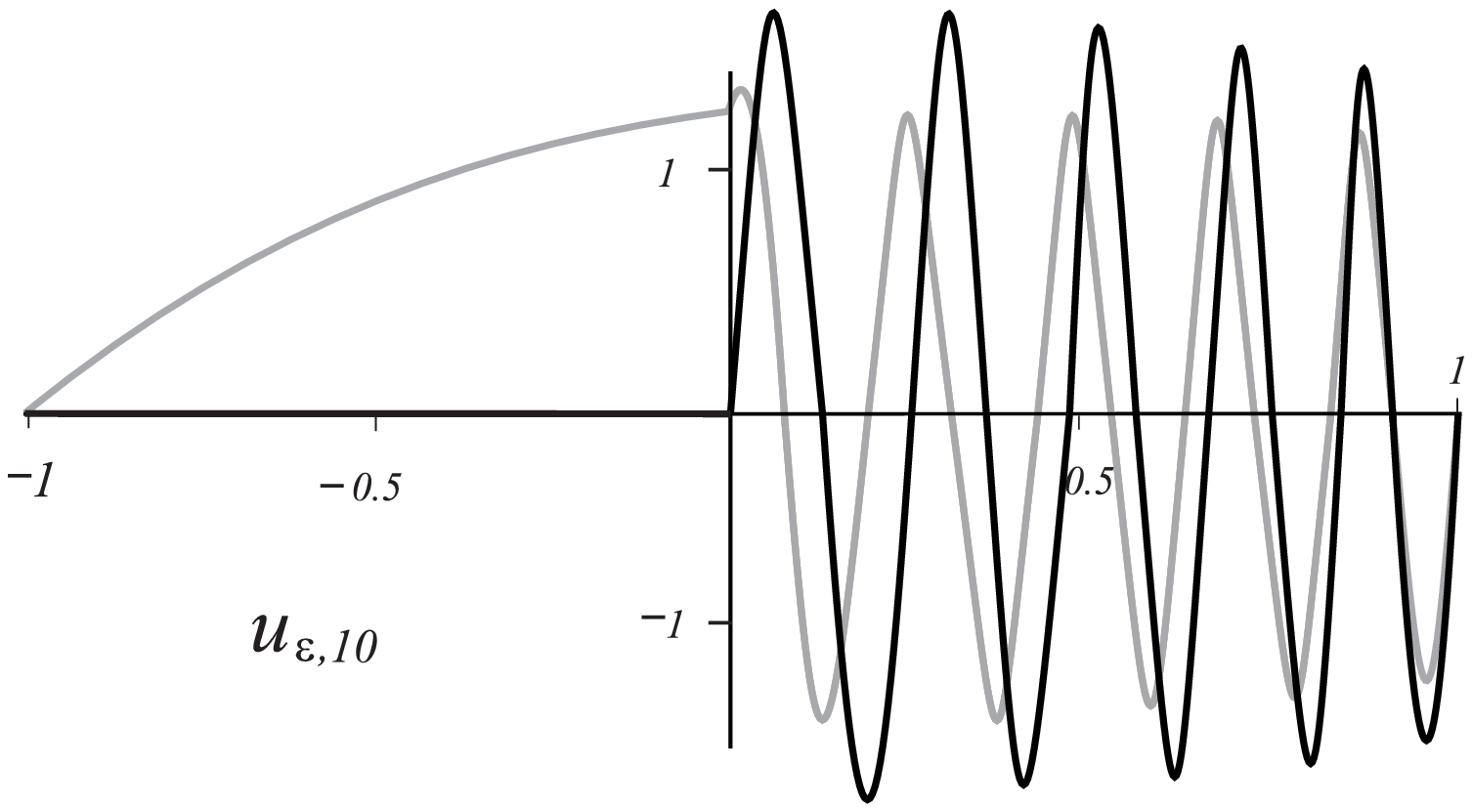}%
\hglue40pt\raise6.6pt\hbox{%
\includegraphics[width=0.4\textwidth]{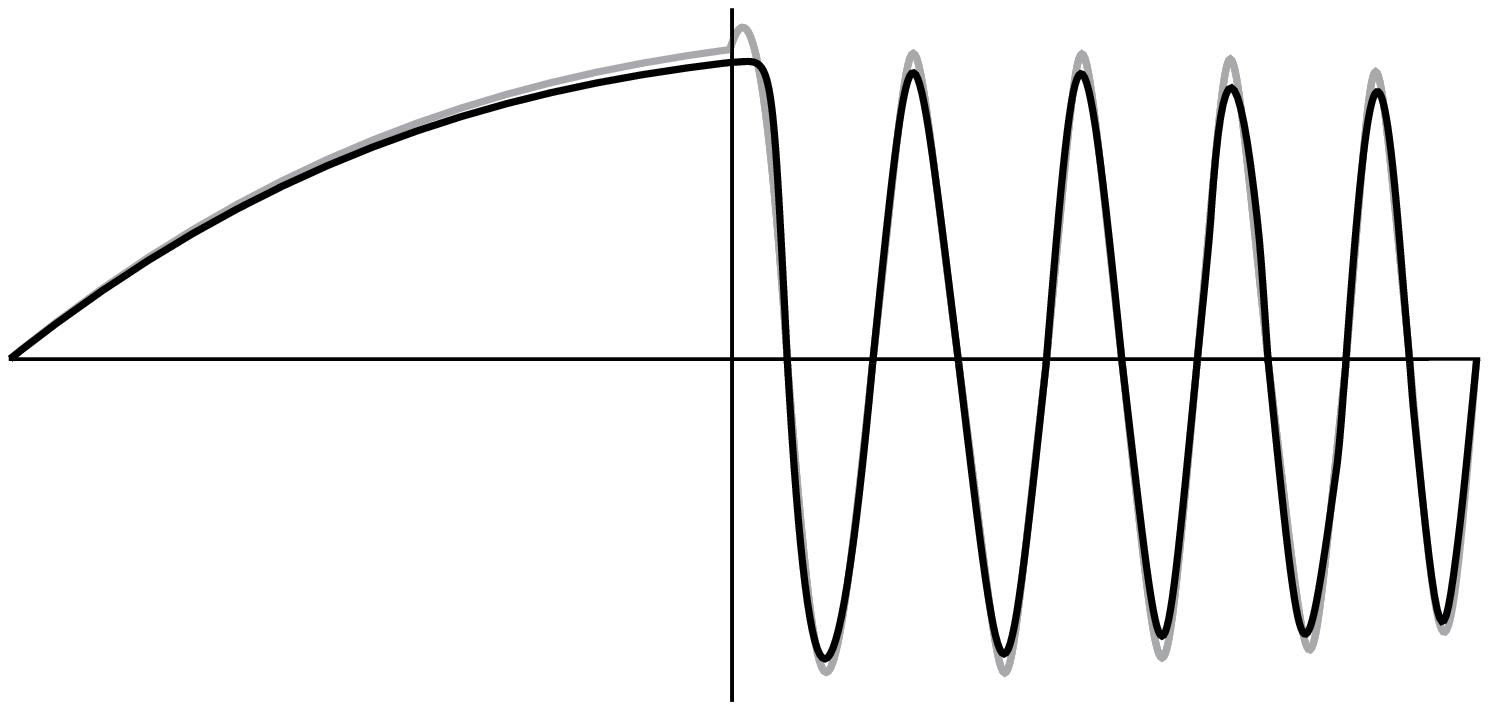}}
\\[-5pt]
\vglue12pt
   \includegraphics[width=0.4\textwidth]{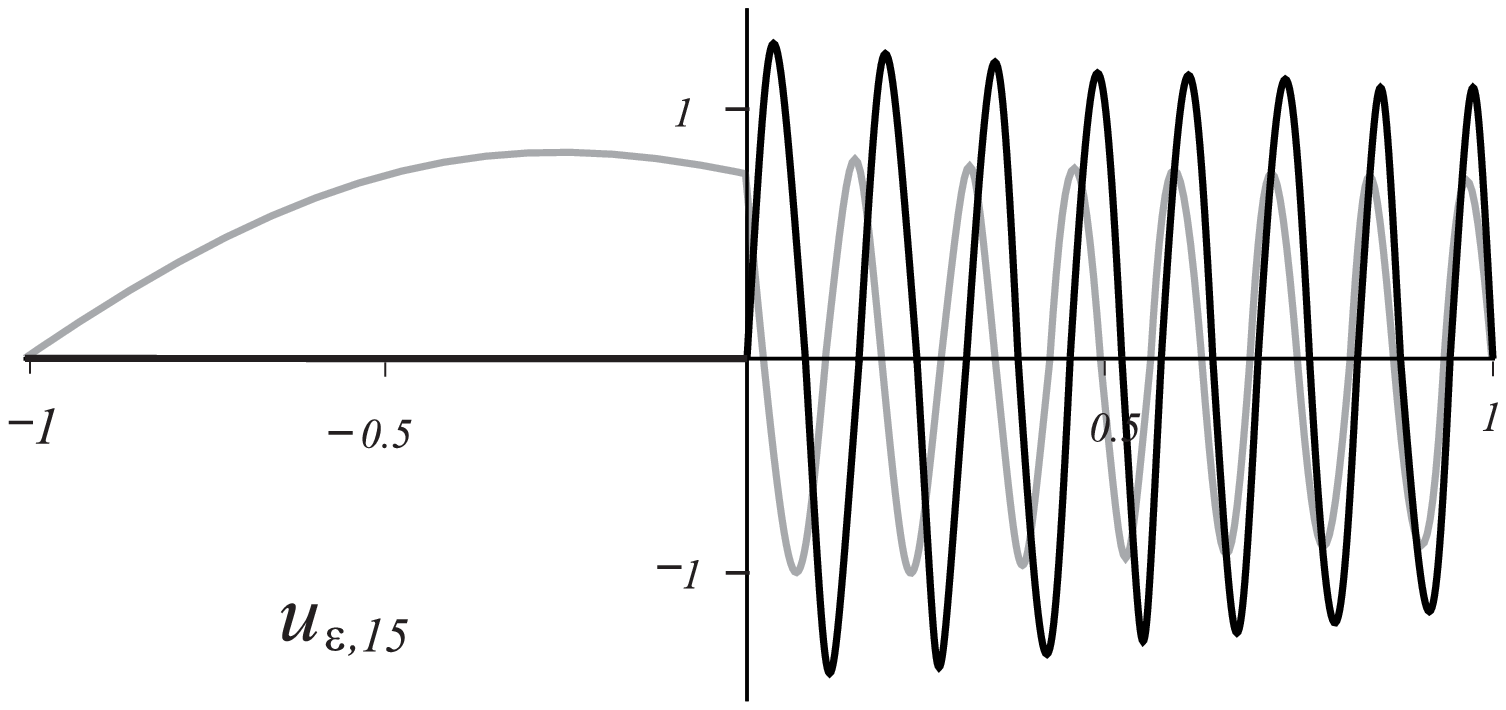}%
   \hglue40pt \includegraphics[width=0.4\textwidth]{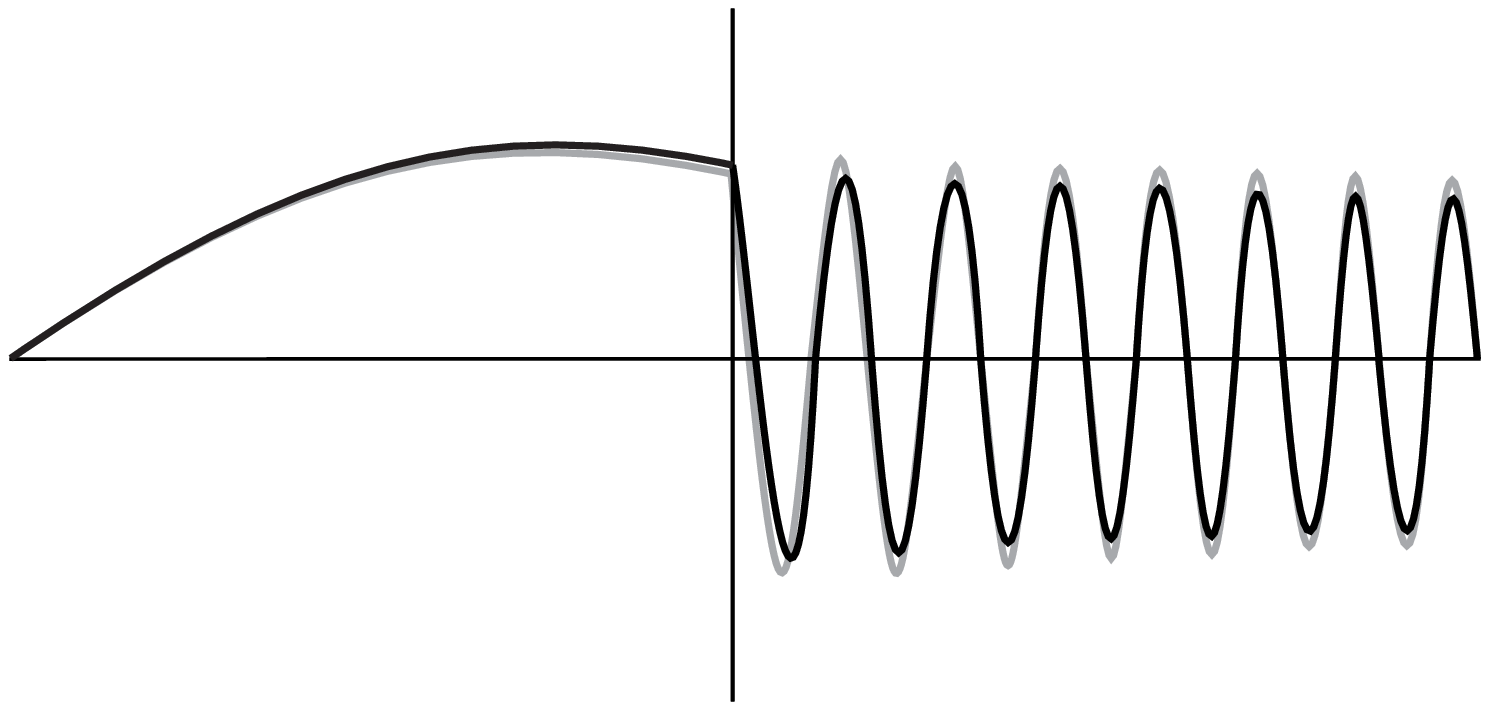}
\caption{\it A comparison of the low and high frequency approximations (black plots) with the 
eigenfunctions (grey plots) for $u_{\eps,5}$, $u_{\eps,10}$ and $u_{\eps,15}$ (from top to bottom). }
\end{center}
\end{figure}

Note that the method that is applied for the approximations
in $1$--dimensional case
is also applicable
in a multidimensional situation.
Nevertheless, the justification of it
requires another technique, which
is out of scope of this paper.


\end{document}